\documentclass[]{article}

\usepackage{amsmath,amssymb,times}
\usepackage{graphicx}
\usepackage{epstopdf}
\usepackage{color}

\def\integers{\ensuremath\mathbb{Z}}
\let\Z=\integers
\newcommand{\zn}{\mathbb Z_n}
\def\Q{\ensuremath\mathbb{Q}}
\newcommand{\C}{\ensuremath \mathbb C}

\definecolor{blue}{rgb}{0,0,1}   \definecolor{red}{rgb}{1,0,0}

\newcommand{\Diag}{\text{ Diag }}
\newcommand{\lcm}{\text{lcm}}

\def\matplural{matrices}

\let\mc=\mathcal

\newtheorem{Def}{Definition}
\newtheorem{example}{Example}
\newtheorem{Prop}{Proposition}
\newtheorem{Thm}{Theorem}
\newtheorem{Lemma}{Lemma}
\newtheorem{remark}{Remark}

\title{On the group of rational spectral units with finite order}

\begin{document}
\author{Emmanuel Amiot}\thanks{manu.amiot@free.fr}
\date{5/21/2009, CPGE, Perpignan, France }

\maketitle

\abstract{ The problem of phase retrieval is a difficult one which remains far from solved. Two homometric sets are always connected by way of a convolution product by some {\em spectral unit}, though not necessarily in a unique way. 
Here we elucidate one small aspect, the subgroup of  spectral units with finite order. Its elements are completely characterized by relations between their eigenvalues. This sheds some light on the beltway problem.}\medskip

\medskip

\section{Introduction}
In this paper we fix some integer $n$ and consider vectors in $\C^n$ as maps from the cyclic group $\zn$ to $\C$:
$a = (a(0), \dots a(n-1))  = (a_0 \dots a_{n-1})$. Usually, such a vector can be the characteristic function of a subset of $\zn$, e.g. in $\Z_{12}$ we put $a = (1,0,0,0,1,0,0,1,0,0,0,0)$ for subset $\{0,4,7\}$.
Two vectors $a,b$  are {\em homometric} if equivalently
\begin{itemize}
   \item
   They have the same interval distributions (i.e., histogram).
   \item
   Their Fourier coefficients have the same magnitudes.
\end{itemize} 
This problem originated with diffraction patterns, as illumination of such patterns is essentially
proportional to the amplitude of a Fourier transform. If we know that $a$ is homometric to some given set $b$ then finding $a$ exactly is known as the {\em phase retrieval problem}, because already the amplitude of $a$'s Fourier coefficients is known. The problem is that $a$ is not necessarily congruent to $b$, as seen from the example:
\begin{example}
   Consider the characteristic functions of $\{0,1,4,6\}$ and $\{0,1,3,7\}$, subsets of $\Z_{12}$ i.e.
$a=(1,1,0,0,1,0,1,0,0,0,0,0), b=(1,1,0,1,0,0,0,1,0,0,0,0)$. The sets are not congruent but are homometric:
for instance the interval distributions are the same, e.g. (4,1,1,1,1,1,2,1,1,1,1,1). Alternatively the Fourier coefficients have magnitudes $(4,\sqrt 2, 2, \sqrt 2, 2, \sqrt 2, 2, \sqrt 2, 2, \sqrt 2, 2, \sqrt 2)$.
\end{example}

In a paper with William Sethares \cite{AmiotSethares} about musical scales or rhythms expressed in terms one of another, we introduced {\em scale matrices}, which are circulating matrices with a given first column.

\begin{Def}
 $M\in\mathcal M_n(\C)$ is a circulating matrix iff its general coefficient is $m_{i,j} = \alpha_{i-j}$
 (where $i-j$ is computed modulo $n$).
 \end{Def}
Hence any circulating matrix $M$ can be written $M : \alpha_0 I_n + \alpha_1 J + \dots \alpha_{n-1} J^{n-1}$
where the vector $\alpha = (\alpha_0, \alpha_1, \dots \alpha_{n-1} )$ is given and
 $ J = \left( \begin{smallmatrix}
       0 & 0 & \dots & &   1 \\
       1 & 0 & 0 \dots & &   0\\
       0 &1 & \dots & &  0\\
       \vdots & & \ddots& &\vdots\\
       0 &  & \dots & 1 & 0
    \end{smallmatrix} \right)
$.

The scale matrix associated with the scale vector $a=(a_0, a_1, \dots a_{{n-1}})$ is precisely
$$
  M(a) = a_0 I_n + a_1 J + \dots a_{n-1} J^{n-1} = 
  \begin{pmatrix}
       a_0 & a_{n-1} & \dots & &   a_1 \\
       a_1 & a_0 &  \dots & &   a_{n-1}\\
       a_2 &a_1 & \dots & &  \vdots\\
       \vdots & & \ddots& &\vdots\\
       a_{n-1} &  & \dots & a_1 & a_0
    \end{pmatrix}.
$$
All these matrices form a commutative algebra $\mc C_n$. 

The matrix product $M(a)\times M(b)$ corresponds with the convolution product of $a$ and $b$, as seen from the first column of the product matrix:
$$
   M(a)\times M(b) = M(a\star\, b),
   a\star b = (\sum_{k=0}^{n-1} a_k b_{-k}, \sum_{k=0}^{n-1} a_k b_{1-k}, \dots )
$$
where all indexes are computed modulo $n$.

The connection with Fourier coefficients is straightforward:
\begin{Prop}
  We define the Fourier matrix $\Omega$ by its $(j,k)^{th}$ term, equal to $\dfrac1{\sqrt n} e^{2 i j k \pi/n}$.
  Notice that ${}^t \Omega = \Omega$ and $\Omega^{-1} = \overline\Omega$. Then
  for any $a\in\C^n, \Omega^{-1} M(a) \Omega$ is equal to the diagonal matrix $F_a$
   with coefficients 
   $$\mc F_a(0) = \sum_{k=0}^{n-1} a_k,\dots \mc F_a(j) = \sum_{k=0}^{n-1} a_k e^{2 i j k \pi/n}, \dots
   \mc F_a(n-1) = \sum_{k=0}^{n-1} a_k e^{2 i (n-1) k \pi/n} = \sum_{k=0}^{n-1} a_k e^{- 2 i k \pi/n}
   $$ 
   i.e. the Fourier coefficients of $a$.
\end{Prop}
We retrieve from there the famous relation $\mc F(a\star b) = \mc F(a) \times \mc F(b)$.\footnote{ This explains why the equality of the interval contents is equivalent to the equality of amplitudes of Fourier coefficients, as the interval content is computed from $a\star (-b)$.}
In more abstract terms, three algebras are isomorphic:
\begin{enumerate}
 \item
   the algebra of columns $\mathbb C^n$, identified with the maps from $\mathbb Z_n$ to $\mathbb C$, and with the convolution product; 
   \item
   the algebra of circulating \matplural, and
   \item 
    the algebra of diagonal \matplural.
\end{enumerate}
Homometry is then readable on the matrix with the magnitudes of the Fourier coefficients, which is (squared)
$$
   \Omega^{-1} M(a) \Omega. {}^t\overline{\Omega^{-1} M(a) \Omega}
   = \Omega^{-1} M(a){}^t \overline{M(a)} \Omega
$$
For real valued vectors $a$ we have of course $\overline{M(a)} =M(a)$.
In other words, $a$ and $b$ are homometric iff
$a\star a^* = b\star b^*$ i.e. $ M(a){}^t \overline{M(a)} =  M(b){}^t \overline{M(b)}$.
(NB: $a^*$ is derived from $a$ by retrogradation and conjugation)\medskip

Another way to look at it is to consider a diagonal matrix $D$ such that
$$
    F_b = \Omega^{-1} M(b) \Omega = D\,  \Omega^{-1} M(a) \Omega = D F_a
$$
where the diagonal elements of $D$ are of amplitude 1 for $a, b$ to be homometric. We get in return a characterization:
\begin{Prop}
   $a,b$ are homometric if there exists a matrix $U$ such that
   \begin{itemize}
   \item
   $ {}^t \overline U = U^{-1}$
   \item 
   $U\in\mc C_n$
   \item
   $B = U A$.
   \end{itemize}
\end{Prop}
Equivalently, we have $b  = u\star a$ with $u\star u^* = (1,0,0,\dots)$ (the more usual definition, see \cite{Rosenblatt}).
Such elements satisfying $u\star u^* = (1,0,0 \dots) = \delta_0$, i.e. homometric to $ \delta_0$, are called {\em spectral units}.

The difficult issue there is: where do the coefficients of $u$ (or $U$) lie ? Abstractly, the group of all possible $U$'s is a torus, as each eigenvalue can be any complex number on the unit circle (the eigenvectors are always the same, since they are shared by all circulating matrixes).
In practice (for instance in music theory), we look for integer, or at worst rational, coefficients.
Rosenblatt proved (his Thm. 3.1, \cite{Rosenblatt}) that  
\begin{Prop}
  For any given subfield $K$ of $\C$, closed under conjugation, 
  for $a, b\in K^n$, $a$ and $b$ are homometric iff there exists a spectral unit $u\in K^n$ satisfying
  $b = u\star a$.
\end{Prop}
For instance, $u=(0,1,0,0\dots 0)$ gives rise to a permutation circular matrix $U$. As a spectral unit it realizes the translation of a subset by 1. It is worthy of note that symmetries of the cyclic group $\zn$, mapping the characteristic function of some subset $S$ to the characteristic function of, say, $-S$, are much more complicated in terms of spectral units. For instance the spectral unit connecting 
$$1_{\{0,3,7\}} = (1,0,0,1,0,0,0,1,0,0,0,0)$$
 and
$$1_{7-\{0,3,7\}} = 1_{\{0,4,7\}} = (1,0,0,0,1,0,0,1,0,0,0,0)$$ (minor and major triads in music theory), is 
$$(7,4,-2,1,7,4,-2,1,-8,4,-2,1)/15.$$

In the context of this paper, it means that $B = U A$ where (taking henceforth $K=\Q$)
\begin{itemize}
   \item
   $U$ is unitary: ${}^t \overline U = U^{-1}$,
   \item
   $U\in \mc C_n$, and
   \item
   $U$ is rational valued: $U\in\mc M_n(\Q)$
\end{itemize}
NB: the result is easy if $A$ (and hence $B$) are non singular, because $U = B A^{-1}$.

We call such \matplural \ {\em spectral units} again.

Unfortunately, as the example above with the major and minor triad shows, some spectral units are not simple -- the one in the example has infinite order. 
This paper enumerates all (rational) spectral units with finite order. Since even in the unit circle considered as a group, both elements with a finite order and rational elements form infinite subgroups, it is notable that these spectral units are finitely many. 
To be specific we will show three theorems. Remember we are dealing with rational spectral units with finite order only.
 \begin{Thm} \label{order}
    If $U$ is a spectral unit with finite order and $n$ is even, then all its eigenvalues are $n^{th}$ roots of unity.
    If $n$ is odd, then the eigenvalues are either $n^{th}$ roots of unity or their opposites (i.e. $2n^{th}$ root of unity).
 \end{Thm}
This stems from a more precise condition
 \begin{Thm} \label{condition}
    If $U$ is a spectral unit with finite order and $n$ is even, then for all $k$ coprime with $n$ and 
    any Fourier coefficient  (= eigenvalue of $U$)$\xi_j, j\neq 0$, one has $\xi_{k j} =\xi_j^k$. 
    For $j=0$ we have $\xi_0 = \pm 1$.
    
   This condition stands even when $n$ is odd, except when $\xi_j$ is a $e^{(2p_j +1)i\pi/n}$ and $k$ even, when we have
   $\xi_{k j} = -\xi_j^k$.
 \end{Thm}
 For instance for $k=-1$ this gives the condition that the last Fourier coefficients must be the conjugates of the first ones (thus ensuring that $U$ is real valued).

These conditions enable to specify exactly what are the possible values of the Fourier coefficients $\xi$, and hence to retrieve all possible spectral units with finite order (by inverse Fourier transform). We will see some examples in the next section. A general result can be stated which gives the precise cardinality of the group of such spectral units, as it gives directly the decomposition of this (abelian) group into a product of cyclic groups.

 \begin{Thm} \label{structure}
    Any spectral unit with finite order is defined by the values of the subset $\{\xi_j, j\mid n\}$ of its eigenvalues.
    The possibilities are listed {\em infra}:
    \begin{itemize}
      \item
      $\xi_0=\pm 1$;
      \item
       When $n$ is odd, for all $j\mid n$, 
       $\xi_j$ OR $ -\xi_j$ is any power of $e^{2i j\pi/n}$.
       \item
       When $n$ is even, $\xi_j$ is any power of $e^{2i j\pi/n}$ if $n/j$ is even, or any power of $e^{i j\pi/n}$ if $n/j$ is odd. 
    \end{itemize}
 \end{Thm}
 
 Unfortunately, many spectral units are not of finite order (even in simple cases), like
 $(0, 1/2, 0, 1/2, 1/2, -1/2, 0)$. Still the results above may perhaps enable to compute all spectral units with, say, small denominators, which occur in practice for homometric subsets of $\zn$.
 
 \section{Examples}
 
 \begin{example}
    Let us elucidate the group of units when $n=12$. Let $u$ be a spectral unit with finite order, and $\xi_0, \dots \xi_{11}$ its Fourier coefficients, i.e. eigenvalues of the associated matrix $U$. From theorem 2 above,
    the relation $\xi_j^k = \xi_{j k}$ is satisfied for all four values of $k = 1,5,7,11$.
    \begin{itemize}
      \item 
      There are no conditions on $\xi_1$ which is any $12^{th}$ root of unity; its value specifies $\xi_5 = \xi^5$ and similarly $\xi_7, \xi_{11}$.
      \item
      $\xi_2$ must be a power of $e^{2\times 2 i \pi/12} =e^{i\pi/3}$: without making use of Thm. 3, we can see that 
      if we set $\xi_2 = e^{2i\pi\alpha/12}$ where $\alpha\in\integers_{12}$, then
      $\xi_2^5 = \xi_{10} = \xi_{-2} =\xi_2^{-1}$, hence $5\alpha = -\alpha \mod 12$, i.e. $\alpha$ is even.
      This determines also $\xi_{10} = \overline \xi_2$.
      \item
      Similarly $\xi_3$ is a power of $i = e^{i\pi/2}$: if we set again $\xi_3 = e^{2i\pi\alpha/12}$ we get
      $$
         \xi_3^7 = \xi_9 = \xi_{-3} = \xi_3^{-1} \quad\text{i.e.}\quad
              7\alpha = -\alpha \mod 12
      $$
      meaning that $8\alpha=0$ i.e. 3 divides $\alpha$. We have $\xi_9 = \overline \xi_3$.
      \item 
      $\xi_4$ gives the slightly more complicated case:
      it must be a power of $e^{i\pi/3}$, just like $\xi_2$, and not a power of 
       $e^{{\mathbf 2} i\pi/3}  =  e^{2\times 4 i\pi/12}$: if $\xi_4 = e^{2i\pi\alpha/12}$,
       $$
         \xi_4^5 = \xi_8 = \xi_{-4} = \xi_4^{-1} \quad\text{i.e.}\quad
              5\alpha = -\alpha \mod 12
      $$
      Here also, we find that $\xi_8 = \overline \xi_4$.
      \item 
      $\xi_0 \pm 1$ and $\xi_6 = \xi_{-6} = \xi_6^{-1}$ is a $12/6^{th}$ root of 1, i.e. $\xi_6 = \pm 1$.
    \end{itemize}
    To conclude: $\xi_1$ is any $12^{th}$ root of unity, while $\xi_2, \xi_3, \xi_4$ are limited to subgroups, $\xi_0$ and $\xi_6 = \pm 1$. The structure of the group is then
    $\Z_{12} \times (\Z_6)^2\times \Z_4 \times (\Z_2)^2$, with 6,912 elements.
\end{example}

\begin{example}
  Take $n=7$. Then apart from $\xi_0 = \pm 1$, all Fourier coefficients are powers of $\xi_1$: 
  \begin{itemize}
    \item
    If $\xi_1^7 = 1$ then $\xi_2 = \xi_1^2$ and more generally $\xi_k = \xi_1^k, k=2\dots 6$.
    \item
    If $\xi_1^7\neq 1 = \xi_1^{14}$ then we play the same thing with $-\xi_1$:
    $\xi_2 = -\xi_1^2, \xi_4 = -\xi_1^4$ and $\xi_6 = \xi^{-1} (=\xi_2^3), \xi_5 = \xi_2^{-1} = \xi_1^5, \xi_3 = \xi_4^{-1}$.
    
    In the exponential notation, we get more simply
    $$
      (\xi_1, \xi_2, \xi_3, \xi_4, \xi_5, \xi_6) = (-e^{\frac{(2k+1)i\pi}7}, -e^{\frac{2(2k+1)i\pi}7}, \dots -e^{\frac{6(2k+1)i\pi}7})
    $$
     When $\xi_0 = +1$ this gives spectral units of the form 
     $$u = (2/7, -5/7, 2/7, 2/7, 2/7, 2/7, 2/7)$$
     or a circular permutation thereof; if $\xi_0 = -1$ then we obtain opposites of circular permutation matrices, 
     e.g. $u = (0, -1, 0, 0, 0, 0, 0)$.
     These spectral units are exactly the opposite of the ones obtained in the first case, when $\xi_1^7 = 1$.
  \end{itemize}
  The whole group has only 28 elements (its quotient by $\pm 1$ is the cyclic group of $14^{th}$ roots of unity). 
\end{example}

\begin{example}
  We end these examples with an odd composite number: $n=15$. As usual $\xi_0 = \pm 1$.
 \begin{itemize}
   \item
   Coefficients $\xi_k, k\in\Z_{15}^*$, are deduced from $\xi_1 = \xi$:
   \begin{itemize}
     \item
     If $\xi^{15}=1$ (e.g. $\xi = e^{8i\pi/15}$) then 
    \begin{multline*}
        \xi_2 = \xi^2, \xi_4 = \xi^4, \xi_7 = \xi^7,  \xi_8 = \xi^8, \\
      \xi_{11} = \xi^{11}, \xi_{13} = \xi^{13},
        \xi_{14} = \xi^{14} = \xi^{-1} = \overline{\xi}
    \end{multline*}
     e.g. $\xi_k = e^{8ki\pi/15}$ for $k = 2,4,7,8,11,13,14$.
     \item
     If $\xi^{15} = -1$ instead, then the list is the opposite of the one found in the first case, e.g.
     $\xi_k =- e^{8ki\pi/15}$ for $k = 2,4,7,8,11,13,14$.
     \end{itemize}
     \item
     Similarly, $\xi_6, \xi_9, \xi_{12}$ are powers of $\xi_3$ (or their opposites), which must be a $5^{th}$ root of 1 (or -1):
     for example, 
     $$
       \xi_3 = e^{4i\pi/5}, \xi_6 = e^{8i\pi/5}, \xi_{12} =  e^{16i\pi/5} =  e^{6i\pi/5}, 
       \xi_9 = \xi_3^8 = e^{32i\pi/5}= e^{2i\pi/5} 
     $$
     or (when $\xi_3$ is a fifth root of -1) the opposite:
     $$
       \xi_3 = -e^{4i\pi/5}, \xi_6 = -e^{8i\pi/5}, \xi_{12} =  -  e^{6i\pi/5}, 
       \xi_9 = - e^{2i\pi/5} 
     $$
     \item
     The same goes for $\xi_5 = \pm e^{2ik\pi/3}$ and $\xi_{10} = \xi_5^{-1}$.
  \end{itemize}
  More concisely: $\xi_1$ is any $30^{th}$ root of unity, $\xi_3$ any $10^{th}$ root, $\xi_5$ any $6^{th}$ root.
  The group has 
  $2 \times {30} \times {10} \times {6} = 3,600$ elements.
\end{example}

\begin{example}
    Say $a, b $ are homometric. If matrices $A, B$ are non singular, then there is one and only one spectral unit $u$,
    with matrix $U = B A^{-1}$ such that $B = U A$ (or equivalently $b = a \star u$).
    But sometimes eigenvalues of $A, B$ are nil. \cite{Rosenblatt} chooses to give 1 as eigenvalues for the corresponding eigenvectors, but other choices are possible, if we respect the conditions found in the theorems above.
    
An example issued from music theory is $a = (1,0,1,1,0,1,0,1,0,1,0,1), b=(1,0,1,0,1,0,1,0,1,1,0,1)$, two melodic minor scales. The Fourier coefficients with indexes 2 and 10 are nil.
\begin{itemize}
 \item
 Using Rosenblatt's choice, we take $\xi_2 = \mc F_u(2) = \mc F_u(10) = \xi_{10} = 1$ (the other Fourier coefficients are determined by
 $\mc F_u(k) = \mc F_b(k) /\mc F_a(k) $). This yields 
 $u = (0, 0, 0, 0, 0, 0, 0, 0, 0, 1, 0, 0)$. Musically this means that A minor is transposed from C minor by a minor third.
 \item
   We know from Thm. \ref{condition} that $\xi_2$ is some power of $e^{i\pi/3}$, $\xi_{10}$ being its conjugate or inverse. This yields five other possible units, e.g. 
\begin{multline*}
     u = \dfrac14 (1, 0, -1, -1, 0, 1, 1, 0, -1, 3, 0, 1) \text{  or  }\\
     u = \dfrac1{12}(1, 2, 1, -1, -2, -1, 1, 2, 1, 11, -2, -1)\text{  or  }\\
     u = \dfrac16 (2, 1, -1, -2, -1, 1, 2, 1, -1, 4, -1, 1)\text{  or  }\\
     u = \dfrac1{12}  ({1, -1, -2, -1, 1, 2, 1, -1, -2, 11, 1, 2})\text{  or  }\\
      u = \dfrac14 ({1, 1, 0, -1, -1, 0, 1, 1, 0, 3, -1, 0}).
  \end{multline*}
\end{itemize}
  In a way this can be interpreted as other, hidden symmetries between those two musical scales.  
\end{example}

\section{Proofs}

\subsection{Proof of Thm. \ref{order}}
Throughout,  $U$ is a circulating matrix which is unitary ($U^{-1} = {}^t \overline{U}$), is of finite order:
$U^m = I_n$ for some $m$,
and has rational elements. Hence its eigenvalues are of magnitude 1 (they are $m^{th}$ roots of unity), and as discussed above $U$ diagonalizes into
$\Diag(\xi_0, \xi_1, \dots \xi_{n-1})$ where the eigenvalues $\xi_j$ are also the Fourier coefficients of the first column of $U$, seen as a map from $\zn$ to $\C$.

We prove an alternative form of Thm. \ref{order}:
\begin{Prop}
    All eigenvalues of $U$ are $n^{th}$ roots of unity for even $n$, and $2n^{th}$ roots of unity for odd $n$.
\end{Prop}
As we assume $U$ has finite order, all these eigenvalues  are roots of unity.
Moreover, as $U$ is a polynomial $U = P(J), P\in\Q[X]$ in the matrix $J$, whose eigenvalues are the $n^{th}$ roots of unity, the eigenvalues of $U$ are polynomials in these roots: $\xi_k =P(e^{2ik\pi/n})$ , hence lie in the field $\Q[e^{2i\pi/n}] = \Q_n$, often called {\em cyclotomic field}.
We need the following
\begin{Lemma}\label{order}
   Let $\xi$ be a $m^{th}$ root of unity belonging to the cyclotomic field $\Q_n$.
   
   Then $\begin{cases}
      \xi^n = 1 & \text{when $n$ is even,}\\
      \xi^{2n} = 1 & \text{when $n$ is odd.}
   \end{cases}$
\end{Lemma}
In other words, if $\Q_m\subset \Q_n$ then $m$ is at most $n$ or $2n$, according to whether $n$ is even or odd.\footnote{ For instance $\Q_3 = \Q_6$.}

Let $\xi$ be such a number (like any eigenvalue of $U$). Let $m$ be the order of $\xi$, i.e.  the smallest integer satisfying $\xi^m=1$; we know that $\xi$ generates  $\Q_m$.
As  $\xi\in\Q[e^{2i\pi/n}]$ too, $\Q_m \subset \Q_n$. This does not preclude $m>n$. We need the following
\begin{Lemma}
   The multiplicative group of elements of finite order in $\Q_n$ is cyclic.\footnote{ It is perhaps not obvious that this group is finite, and indeed the group of elements of $\Q_n$ with length one is not; this holds because for large $m$ the dimension of $\Q_m/\Q$ exceeds that of $\Q_n/\Q$, equal to $\varphi(n)$, cf. below.}
\end{Lemma}
This is because given two elements $\xi, \xi'$ with orders $m, m'$ it is possible to construct an element of order $\lcm (m, m')$. In other words, the roots of unity in $\Q_n$ have a maximum order, which is the lcm of all possible orders.

 Let us call again $m$ this maximal value, to prove Lemma 1 we need to prove that $m=n$ or $2n$.
 Now, any element $\xi$ of $\Q_n$ which is a root of unity must satisfy $\xi^m = 1$. 

This is true in particular when $\xi$ is the primitive $n^{th}$ root $e^{2i\pi/n}$;
hence $m$ is a multiple of $n$, which entails $\Q_n \subset \Q_m$.
Finally $\Q_n = \Q_m$.

Letting $\varphi(n) = \dim (\Q_n / \Q)$ stand for Euler's totient function, 
$
  \begin{cases} n \mid m \text{ and}&  \\\ \varphi(n) = \varphi(m)  & \ \end{cases}.
$

As $\varphi(n) = n\prod\limits_{p\mid n; p \text{ prime}} \left(1 - \dfrac1p\right)$, the only possibility is that $m = \begin{cases} n & \text{ for $n$ even} \\ 2n & \text{ for $n$ odd} \end{cases}$.
This proves the Lemma \ref{order}, and hence the Proposition: all eigenvalues of $U$ are $n$ or $2n^{th}$ roots of unity.
Let us clarify the case of odd $n$: $e^{i\pi/n} = - {\bigl(e^{2i\pi/n}\bigr)}^{\frac{n+1}2}$
and we do have $\Q_n = \Q_{2n}$. So we can rephrase the Lemma: in the odd case, 
$\xi^n = \pm 1$.\medskip

\begin{remark}
At this point, $U$ could be constructed as a polynomial in the elementary circulating matrix $J$
   (as all other circulating \matplural) $U = P(J)$, where $P$ is the interpolating polynomial that sends
   the Fourier coefficients of $J$, i.e. the $e^{2ik\pi/n}$, to the Fourier coefficients chosen for $U$.
   Such a construction is easy, with the Lagrange polynomials associated with the $e^{2ik\pi/n}$, as
   $P$ is a linear combination of these polynomials with coefficients that are precisely 
   the Fourier coefficients of the desired $u$.
\end{remark}

\subsection{Proof of Thm. \ref{condition}}

The possibilities of mapping the $n^{th}$ roots of 1 to $m^{th}$ roots of 1 can be somewhat reduced by noticing that
$U$ is a rational polynomial\footnote{ The coefficients of $P$ can be read on the first column of $U$.} in $J$, and such a polynomial is stable under all field automorphisms of $\Q_n$ if we use the following characterization from Galois theory:
\begin{Lemma}\label{fixed}
     Any object (number, vector, polynomial, matrix) with coefficients in $\Q_n$ is rational-valued iff it is invariant under all Galois automorphisms of the cyclotomic extension $\Q_n $ over $\Q$.
\end{Lemma}
We mention without proof either the structure of the Galois group. 
These two results can be found in any textbook on Galois theory.
\begin{Lemma}
   Any field automorphism of the cyclotomic extension $\Q_n$ over $\Q$ is defined by 
   $\Phi_k(e^{2i\pi/n}) = e^{2ik\pi/n}$ for some $k\in\zn^*$, the group of invertible elements of the ring $\zn$, e.g. for any integer $k$ coprime with $n$.
\end{Lemma}
This is enough to define $\Phi_k(x)$ for any $x\in\Q_n$, as $x = \sum a_j e^{2i j\pi/n}$ with rational $a_j$'s, and hence
$\Phi_k(x) = \sum a_j e^{2i j k\pi/n}$.
For instance when $n=12$, the four different automorphisms $\Phi_k$ are defined by the possible images of 
$e^{ 2i \pi /12} = e^{ i \pi /6}$, namely
 $e^{ i k\pi /6}, k\in\{1,5,7,11\}$. Their group is isomorphic with the multiplicative group $\Z_{12}^* = \{1,5,7,11\}$.
 
 If $\Phi_k$ is such an automorphism, notice that $\Phi_k(\xi) = \xi^k$ for any $n^{th}
$ root $\xi$  of unity (with one exception: $\Phi_k( - 1) = -1 \ \forall k\in\zn^*$).
If $n$ is odd and $\xi$ is a $2n^{th}$ root but not a $n^{th}$, then $-\xi$ is a $n^{th}$ root, and hence
$$\Phi_k(\xi) = -(-\xi)^k = \begin{cases}
\xi_k & \text{for odd $k$'s} \\ -\xi_k & \text{for even $k$'s}\end{cases}.$$
For instance  $\Phi_2(\xi)  = -\xi^2$ for such $\xi$.
\medskip

 So from Lemma \ref{fixed}, we will have that $U\in \mathcal M_n(\Q)$ if $U$ is invariant under all the $\Phi_k, k\in\zn^*$.\medskip

Now we can prove Thm \ref{condition}, beginning with the case of $n$ even.
\smallskip

 Consider the eigenvector $X_j = (1, e^{2i j \pi /n},  e^{2i 2j \pi /n}, \dots  e^{2i j(n-1) \pi /n})'$ for the eigenvalue $\xi_j$ of $U$ (for matrix $J$, the eigenvalue is of course $ e^{2i j \pi /n}$). The prime signals that we consider $X_j$ as a column.
 We have $\Phi_k(X_j) = X_{j k}$ by direct computation.
 The case $j=0$ is straightforward, as the eigenvector is real valued, so must be the eigenvalue, i.e. $\xi_0 = \pm 1$.
 We exclude this case from now on.
 
 We assume that $\Phi_k(U) = U$ (i.e. that $U$ is rational valued).
 Applying the field automorphism $\Phi_k$ to the equation 
 $$ U X_j = \xi_j X_j \quad\text{yields}
 \quad \Phi_k(U) \Phi_k(X_j)  = U X_{k j} = \xi_{k j} X_{k j}
   =\Phi_k(\xi_j) \Phi_k(X_j) = \xi_j^k X_{k j}
 $$
 Hence 
 $$\Phi_k(\xi_j) = \xi_j^k = \xi_{j k} \quad (\sharp)$$ 
 for all $j\neq 0$ and all $k\in\zn^*$.
\medskip
 
 Now for the reciprocal.
 Assume the above equation $(\sharp)$  between the eigenvalues. We choose one Galois automorphism, $\Phi_k$ (for some $k$ coprime with $n$).
Let us apply $\Phi_k(U)$ to any eigenvector $X_j$ of $U$: notice that $X_j = \Phi_k (X_{k^{-1} j})$ where $k^{-1}j$ is computed modulo $n$.
Hence
\begin{multline*}
    \Phi_k(U) X_j = \Phi_k(U X_{k^{-1} j}) = \Phi_k( \xi_{k^{-1} j}X_{k^{-1} j})
    = \Phi_k( \xi_{k^{-1} j}) \Phi_k(X_{k^{-1} j})\\
    = \xi_{k^{-1} j}^k X_{k k^{-1} j} \quad\text{because $\Phi_k$ raises any root of 1 to the $k^{th}$ power} \\
    = \xi_j X_j  \quad\text{by our assumption on the eigenvalues}
\end{multline*} 
We have proved that $\Phi_k(U)$ does the same thing as $U$ on any eigenvector; these make up a basis, hence $\Phi_k(U) = U$, i.e. $U$ is rational valued.
\medskip

It remains to discuss the case of odd $n$. We still get the equation 
$\Phi_k(\xi_j) = \xi_{j k}$ if $U$ is assumed to be invariant under $\Phi_k$.

If $\xi$ is a $n^{th}$ root of unity, the computation is identical.

If $\xi^{2n}= 1$ but  $\xi^n\neq 1$, then $(-\xi)^n=1$ and hence $\Phi_k(\xi) = -\Phi_k(-\xi) = - (-\xi)^k = -\xi^k$ for even $k$ and $\Phi_k(\xi) = \xi^k$ for odd $k$. The computation above still yields 
$\xi_{j k} = \Phi_k(\xi_j) = \xi_j^k$ for odd $k$, and we have also the new case $\xi_{j k} =  -\xi_j^k$ for even $k$.

Say $k=2$, and $\xi_1 = \xi$ with $\xi_1^{2n} = 1\neq \xi_1^n$; 
then $\xi_2 = -\xi^2, \xi_4 = -\xi^4, \dots \xi_{2^m} = -\xi^{2^m}$.
$2$ has a finite order in $\zn^*$, hence for some $m$, $\xi_{2^m} = \xi_ 1$. We get an orbit of $m$ eigenvalues which are all $2n^{th}$ roots of unity, e.g. $\mc O = \{\xi_1, \xi_2, \xi_4, \xi_8 \dots\}$.

Say now that $k = 2^v k', k'$ odd and coprime with $n$. The formula $(\sharp)$ is then valid and yields $\xi_{k} = \xi_{2^v}^{k'}$. So $\xi_k$ is determined when $\mc O$ is known. Notice that $\xi_k$ will never be a $n^{th}$ root (because 2 and $k'$ are coprime with $n$): either all the eigenvalues [with even index] are $n^{th}$ roots, or none (except of course $\xi_0=\pm 1$).

The reciprocal is similar to the even case: it is identical when the eigenvalues are of order $n$ (at most); and if $\xi_1$ has order $2n$ then the values of $\xi_k$ that we have obtained enable to satisfy the relations $ \Phi_k(U) X_j = \xi_j X_j $ for all $j,k$ ¤$k$ coprime with $n$), so that $\Phi_k(U)$ is identical to $U$ i.e. $U$ is rational-valued. This ends the proof.

\subsection{Proof of Thm. \ref{structure}}
We make us of the conditions found in Thm. \ref{condition}.

The whole set of eigenvalues is thus determined if we know $\xi_j$ for a subset of representatives $j$ of all orbits under multiplication by elements of $\zn^*$ (so called {\em associated elements} in the ring $\zn$). We can specify the smallest representatives:
\begin{Lemma}
   Any element $j\in\zn$ is associated with a divisor of $n$, i.e. $\exists k\in\zn^*, k j = \gcd(n, j)$.
\end{Lemma}
(We identify integers and classes modulo $n$ when the distinction is irrelevant).

This stems from the Bezout identity (in $\Z$): for some $k, \ell,  k j + \ell n = \gcd(n, j)$.
After division by $ \gcd(n, j)$ we see that $k$ and $n$ are coprime. But modulo $n,  k j = \gcd(n, j)$, qed.\footnote{ For instance for $n=15$ we have the orbits of equivalent elements 
$$(0), ({\mathbf 1}, 2, 4, 7, 8, 11, 13, 14), ({\mathbf 3}, 6, 9, 12), ({\mathbf 5}, 10)$$ indexed by the divisors $1,3,5$
and of course 0.}

So it is sufficient to specify $\xi_j$ when $j$ is any divisor of $n$.
We will need a last Lemma, which seems interesting in its own right:\footnote{ Though elementary, the result was unknown to the author and does not appear to be readily available in the literature.}
\begin{Lemma}\label{delta}
   The set of differences $\Delta_n = \zn^* - \zn^* = \{ a - b, (a, b)\in(\zn^*)^2\}$
   is $\zn$ when $n$ is odd, $2\zn$ when $n$ is even.
\end{Lemma}
It is straightforward for $n$ prime, for $n$ an odd prime power, and we notice that 
when $n = 2^m$ then $\zn^*$ = odd numbers, so that $\Delta_n=$ even numbers.
The general case now stems from the chinese remainder theorem.
\medskip

We now procede to prove the theorem. Remember that $\xi_0 = \pm 1$.
\subsubsection{When $n$ is even}
In this case all eigenvalues are $n^{th}$ roots of unity. Let $j$ be any strict divisor of $n$.
\begin{itemize}
\item
When $n/j$ is even,  we can produce $k, k'\in\zn^*$ with $k' - k = \dfrac n j\in 2\Z$ from Lemma \ref{delta}.
Hence (noting that $k\equiv k' \mod n$)
$$
   \xi_j^{k + \frac n j} =\xi_j^{k'} = \xi_{j k'} = \xi_{j k} =\xi_j^{k} 
$$
which proves that $\xi_j^{\frac n j} = 1$, i.e.  $\xi_j$ is a power of $e^{2 i j\pi/n}$.
\item
If $n/j$ is odd (meaning that $j$ contains the same power of 2 as $n$), then Lemma \ref{delta} only provides
$k' - k = \dfrac {2n} j$, and the calculation yields
$\xi_j^{2n/j} = 1$, i.e. $\xi_j$ is a power of $e^{i j\pi/n}$, which ends the even case of the theorem.
\end{itemize}
\subsubsection{ When $n$ is odd}
The case when $\xi_j$ is a $n^{th}$ root is identical to the $n$ even (first) case, as from the last Lemma \ref{delta}, we can again produce two elements $k,k'\in\zn^*$ such that $k' - k = \dfrac n j$. and $\xi_j^{k' - k} = 1 = \xi_j^{n/j}$. 
So the spectral unit is determined when we have chosen a $n/j^{th}$ root of unity $\xi_j$  for each divisor $j$ of $n$.\medskip

Now assume that there is an eigenvalue $\xi_j$ which is {\em not} a $n^{th}$ root. 
Then $ -\xi_j$ is  a $n^{th}$ root, and (as $n/j$ is odd) a similar calculation yields for $k'-k = n/j$, with $k,k'\in\zn^*$,
$$
   (-\xi_j)^{k'} = -\xi_{j k'} = -\xi_{j k} =(-\xi_j)^{k} = (-\xi_j)^{k' - \frac n j}
$$
hence $-\xi_j$ is again a $n^{th}$ root of unity.

 This ends the proof of the odd case of Thm. \ref{structure}.

\section*{Acknowledgements}
I thank William Sethares for a very enjoyable common work on scale vectors and matrices which happened to lead to deep questions, among which the one discussed here; Moreno Andreatta, Daniele Ghisi, John Mandereau in Ircam and Mate Matolcsi for a profitable discussion in Paris which encouraged me to explore the matricial traduction of spectral units.

\end{document}